\documentclass[a4paper,english,fontsize=11pt,parskip=half,abstracton]{scrartcl}
\usepackage{babel}
\usepackage[utf8]{inputenc}
\usepackage[T1]{fontenc}
\usepackage[a4paper,left=20mm,right=20mm,top=30mm,bottom=30mm]{geometry}
\usepackage{amsmath}
\usepackage{amsthm}
\usepackage{amssymb}
\usepackage{aliascnt}
\usepackage[bookmarks=false,
            pdftitle={},
            pdfauthor={Benjamin Sambale},
            pdfkeywords={},
            pdfstartview={FitH}]{hyperref}

\newtheorem{Thm}{Theorem} %[section]
\newaliascnt{Lem}{Thm}

\aliascntresetthe{Lem}
\newaliascnt{Prop}{Thm}

\aliascntresetthe{Prop}
\newaliascnt{Cor}{Thm}

\aliascntresetthe{Cor}
\newaliascnt{Con}{Thm}
\newtheorem{Con}[Con]{Conjecture}
\aliascntresetthe{Con}
\theoremstyle{definition}
\newaliascnt{Def}{Thm}

\aliascntresetthe{Def}
\newaliascnt{Ex}{Thm}

\aliascntresetthe{Ex}

\numberwithin{equation}{section}

\allowdisplaybreaks[1]

\renewcommand{\phi}{\varphi}

\newcommand{\Z}{\operatorname{Z}}

\newcommand{\Aut}{\operatorname{Aut}}

\newcommand{\Irr}{\operatorname{Irr}}

\newcommand{\Ker}{\operatorname{Ker}}

\mathchardef\ordinarycolon\mathcode`\:  %defines a nice ":=" 
 \mathcode`\:=\string"8000
 \begingroup \catcode`\:=\active
   \gdef:{\mathrel{\mathop\ordinarycolon}}
 \endgroup

\title{On the projective height zero conjecture}
\author{Benjamin Sambale\footnote{Fachbereich Mathematik, TU Kaiserslautern, 67653 Kaiserslautern, Germany, 
\href{mailto:sambale@mathematik.uni-kl.de}{sambale@mathematik.uni-kl.de}}}
\date{\today}

\begin{document}
\frenchspacing
\maketitle
\begin{abstract}\noindent
Recently, Malle and Navarro put forward a projective version of Brauer's celebrated height zero conjecture on blocks of finite groups. In this short note we show that Brauer's original conjecture implies the projective version.
\end{abstract}

The following is a long-standing conjecture in representation theory of finite groups:

\begin{Con}[Brauer's height zero conjecture~{\cite[Problem~23]{BrauerLectures}}]\label{con1}
Let $B$ be a block of a finite group with defect $D$. Then every irreducible character in $B$ has height $0$ if and only if $D$ is abelian.
\end{Con}

Recently, Malle--Navarro~\cite{MN} proposed the following generalization of \autoref{con1} (the case $Z=1$ yields the original conjecture). An equivalent statement in terms of $\theta$-blocks was given by Rizo~\cite{Rizo}.

\begin{Con}[Malle--Navarro's projective height zero conjecture]\label{con2}
Let $B$ be a $p$-block of a finite group $G$ with defect group $D$. Let $Z$ be a central $p$-subgroup of $G$ and let $\lambda\in\Irr(Z)$. Then every irreducible character in $B$ lying over $\lambda$ has height $0$ if and only if $D/Z$ is abelian and $\lambda$ extends to $D$. 
\end{Con}

In their paper, Malle and Navarro already proved the “if direction” of \autoref{con2} by making use of the solution~\cite{KessarMalle} of the “if direction” of \autoref{con1}. Moreover, they showed that \autoref{con1} implies \autoref{con2} for blocks of maximal defect. 
Generalizing their argument, we prove that \autoref{con1} always implies \autoref{con2}.

\begin{Thm}\label{main}
Suppose that \autoref{con1} holds for all blocks of finite groups. Then \autoref{con2} holds for all blocks of finite groups.
\end{Thm}

Our proof uses the notation from \cite{MN} and the language of fusion systems. Recall that every block $B$ with defect group $D$ induces a (saturated) fusion system $\mathcal{F}$ on $D$ (see \cite[Theorem~IV.3.2]{AKO} for instance). The \emph{focal subgroup} and the \emph{center} of $\mathcal{F}$ are given by
\begin{align*}
\mathfrak{foc}(\mathcal{F})&:=\langle x^{-1}x^f:x\in Q\le D,\ f\in\Aut_{\mathcal{F}}(Q)\rangle\unlhd D,\\
\Z(\mathcal{F})&:=\{x\in D:x\text{ is fixed by every morphism in }\mathcal{F}\}\le \Z(D)
\end{align*}
respectively.

\begin{proof}[Proof of \autoref{main}.]
Let $B$ be as in \autoref{con2}. Since the “if direction” of \autoref{con2} holds, we may assume that $\Irr(B|\lambda)=\Irr_0(B|\lambda)$. By \cite[Theorem~9.4]{Navarro}, the set $\Irr(B|\lambda)$ is not empty and a result of Murai~\cite[Theorem~4.4]{MuraiHZ} implies that $\lambda$ extends to $D$. 
We show by induction on $|G|$ that $D/Z$ is abelian.

Let $K\unlhd G$ be the kernel of $\lambda$. Suppose first that $K\ne 1$. Then $B$ dominates a unique block $\overline{B}$ of $G/K$ with defect group $D/K$ (see \cite[Theorem~9.10]{Navarro}). 
Since the kernel of every $\chi\in\Irr(B|\lambda)$ contains $K$, we have
$\Irr(\overline{B}|\lambda)=\Irr(B|\lambda)=\Irr_0(B|\lambda)=\Irr_0(\overline{B}|\lambda)$. By induction, it follows that $(D/K)/(Z/K)\cong D/Z$ is abelian.

Therefore, we may assume that $\lambda$ is faithful. This implies $D'\cap Z=1$, since $\lambda$ extends to $D$. Let $\mathcal{F}$ be the fusion system of $B$. Then $Z\le\Z(\mathcal{F})$ and it follows from \cite[Lemma~4.3]{DGMP} that $Z\cap\mathfrak{foc}(\mathcal{F})=1$.

Let $\chi\in\Irr(B)$ and $\mu\in\Irr(Z|\chi)$. Since $\mathfrak{foc}(\mathcal{F})\cap Z=1$, there exists an extension $\hat\mu\in\Irr(D)$ of $\mu$ with $\mathfrak{foc}(\mathcal{F})\le\Ker(\hat\mu)$. Similarly, let $\hat\lambda\in\Irr(D)$ be an extension of $\lambda$ with $\mathfrak{foc}(\mathcal{F})\le\Ker(\hat\lambda)$.
By Broué--Puig~\cite[Corollary]{BrouePuigA}, we obtain a character 
\[\psi:=(\hat\lambda\hat\mu^{-1})*\chi\in\Irr(B|\lambda)\] 
(cf. \cite{RobinsonFocal}). By hypothesis, $\psi$ has height $0$ and the same holds for $\chi$, since $\chi(1)=\psi(1)$. Consequently, $\Irr(B)=\Irr_0(B)$ and \autoref{con1} shows that $D$ is abelian and so is $D/Z$.
\end{proof}

\section*{Acknowledgment}
The author is supported by the German Research Foundation (SA \mbox{2864/1-1} and SA \mbox{2864/3-1}).


\begin{thebibliography}{10}

\bibitem{AKO}
M. Aschbacher, R. Kessar and B. Oliver, \textit{Fusion systems in algebra and
  topology}, London Mathematical Society Lecture Note Series, Vol. 391,
  Cambridge University Press, Cambridge, 2011.

\bibitem{BrauerLectures}
R. Brauer, \textit{Representations of finite groups}, in: Lectures on {M}odern
  {M}athematics, {V}ol. {I}, 133--175, Wiley, New York, 1963.

\bibitem{BrouePuigA}
M. Broué and L. Puig, \textit{Characters and local structure in
  {$G$}-algebras}, J. Algebra \textbf{63} (1980), 306--317.

\bibitem{DGMP}
A. D{\'{\i}}az, A. Glesser, N. Mazza and S. Park, \textit{Control of transfer
  and weak closure in fusion systems}, J. Algebra \textbf{323} (2010),
  382--392.

\bibitem{KessarMalle}
R. Kessar and G. Malle, \textit{Quasi-isolated blocks and {B}rauer's height
  zero conjecture}, Ann. of Math. (2) \textbf{178} (2013), 321--384.

\bibitem{MN}
G. Malle and G. Navarro, \textit{The projective height zero conjecture},
  \href{https://arxiv.org/abs/1712.08331v1}{arXiv:1712.08331v1}.

\bibitem{MuraiHZ}
M. Murai, \textit{Block induction, normal subgroups and characters of height
  zero}, Osaka J. Math. \textbf{31} (1994), 9--25.

\bibitem{Navarro}
G. Navarro, \textit{Characters and blocks of finite groups}, London
  Mathematical Society Lecture Note Series, Vol. 250, Cambridge University
  Press, Cambridge, 1998.

\bibitem{Rizo}
N. Rizo, \textit{$p$-blocks relative to a character of a normal subgroup},
  preprint.

\bibitem{RobinsonFocal}
G.~R. Robinson, \textit{On the focal defect group of a block, characters of
  height zero, and lower defect group multiplicities}, J. Algebra \textbf{320}
  (2008), 2624--2628.

\end{thebibliography}
\end{document}